\documentclass[11pt,twoside]{article}

\usepackage{amsfonts}
\usepackage{amsmath}
\usepackage{amssymb}
\usepackage{amsthm}
\usepackage{graphics} 
\usepackage{epsfig} 
\usepackage{latexsym}

\graphicspath{{./figures/}} 

\usepackage[textwidth=16cm,textheight=22cm,top=3cm,bottom=3cm]{geometry}
\usepackage{fancyhdr}

\newtheorem{theorem}{Theorem}

\newtheorem{corollary}{Corollary}

\newtheorem{remark}{Remark}
\theoremstyle{definition}
\newtheorem{definition}{Definition}

\usepackage{color}

\newcommand{\drop}[1]{}
\newcommand{\no}{\noindent}
\newcommand{\fer}[1]{(\ref{#1})}
\newcommand{\qtext}[1]{\quad\text{#1}}

\newcommand{\chil}{\raisebox{0.5mm}{\large $\chi$}}

\newcommand{\cH}{\mathcal{H}}

\newcommand{\eps}{\varepsilon}
\newcommand{\vfi}{\varphi}

\newcommand{\p}{\partial}

\newcommand{\N}{\mathbb{N}}

\newcommand{\R}{\mathbb{R}}

\newcommand{\abs}[1]{| #1 |}

\newcommand{\norm}[1]{\| #1 \|}

\DeclareMathOperator{\argmin}{argmin}

\title{Bifurcation in dynamic problems with seasonal succession
\thanks{Supported by the Spanish MCI Project MCI-21-PID2020-116287GB-I00. 
}}

\author{Gonzalo Galiano and Juli\'an Velasco\thanks{Department of Mathematics, University of Oviedo, Spain 
  ({\tt galiano@uniovi.es, julian@uniovi.es}).}
}

\date{}

\begin{document}

\maketitle

\begin{abstract}
We investigate the bifurcation structure of equilibria in a class of non-autonomous ordinary differential equations governed by a season length parameter, $\tau$, which determines the alternation between growth and decline dynamics. This structure models biological systems exhibiting seasonal variation, such as insect population dynamics or infectious disease transmission. 

Using the Crandall-Rabinowitz bifurcation theorem, we establish the existence of a critical threshold $\tau^*$ at which a bifurcation from the extinction equilibrium occurs. We also explore the emergence of secondary bifurcations from, in general, explicitly unknown non-trivial equilibria which  can only be treated numerically. 

Our results are illustrated with a two-species competitive Lotka–Volterra model for the growth season  and a Malthusian model for the decline season for which primary and secondary bifurcations may be computed analytically, allowing the validation of numerical approximations. 

Our analysis shows how seasonality  drives transitions between extinction of both populations, of only one,  and coexistence of both populations.
  
\vspace{0.25cm}

\no\emph{Keywords: }Seasonal succession; Periodic solutions; Bifurcation; Coexistence; Lotka-Volterra.

\end{abstract}


\section{Introduction}
 We study a class of systems of differential equations ruled by a climatic parameter which determines the equilibria of the system. The climatic parameter splits the temporal line in  seasons where the system behaves in a substancially different manner, e.g. growth and decline.

In its simplest form, the problem is: Find $u:[0,\infty)\to\R^N$ such that $u(0)=u_0$ and
\begin{align}\label{prob}
u'(t)=
\begin{cases}
 f_g(u(t))	& \text{if } t\in G_n,\\
 f_d(u(t))	& \text{if } t\in D_n,
\end{cases}
\end{align}
where the growth and decline seasons are given, respectively, by
\begin{align}\label{def.cndn}
	G_n=(n, n+\tau],\quad D_n = (n+\tau, n+1] \qtext{for } n\in\N, 
\end{align} 
once the period has been normalized to one (e.g. one year). Thus, season duration is captured  by the season length parameter $\tau\in[0,1]$, which is highly influential in the behavior of the solution and which may vary due to climatic change effects \cite{Steinthorsdottir2019,Ruosteenoja2020,Wang2021,Yiqi2023}. 
 In the literature, this problem structure is commonly referred to as \emph{seasonal succession} \cite{Klausmeier2010, Hsu2012, Niu2021, WangZhang2021, Liu2024, Pan2025}, which may be viewed as a sharp limit case of smoother season transitions represented through differentiable time-periodic reaction functions \cite{Namba1986, Peng2012, Lou2017, Browne2020, Chen2023}.

This seasonal structure is a rather common  life framework for many ecological species  
that may grow only during warm seasons while suffering a population decline during cold seasons   \cite{Seeley1985, Almeida2018, Lou2019, Bai2020, Chen2023}.  
Experimental and theoretical work suggests that seasonal variation in resource abundance and differences in seasonal activity behaviors may promote coexistence \cite{Namba1984, Namba1993, Hu1995, Klausmeier2010}, thus not only affecting species growth but also the composition of the ecological niche \cite{DuBowy1988, DeAngelis2009}.

A fundamental question regarding the seasonal succession model \fer{prob} is whether equilibria do exist and, in case they do, what is their dependence on the season length parameter $\tau$. Since $\tau$ measures the growth season duration and $1-\tau$ measures the decline season duration, we may conjecture that a reasonable modeling of these seasons through the functions $f_g$ and $f_d$ should lead to the existence of an extinction equilibrium for $\tau$ small enough.  

However, if $\tau$ is large,  we may also expect the existence of other non-trivial periodic equilibria in which, in general terms,  the solution increases in the growth season and decreases in the decline season (if $u$ is multi-component, this must be taken with care). Thus, we conjecture that for certain types of  functions $f_g$ and $f_d$, a bifurcation may occur at a certain value $\tau_*\in(0,1)$.

To our knowledge, most of the work in the literature has been devoted to the study of particular cases for the growth and decline functions in which explicit solutions for each season may be computed, allowing a detailed study of the global bifurcation of solutions in terms of  $\tau$, among other properties.

For instance, the case of the competitive Lotka-Volterra model for the growth season and the linear Malthusian function for the decline season was introduced and numerically studied by Koch in the early paper \cite{Koch1974}, and analitically elaborated in subsequent works, see \cite{Namba1984, Liu2007, Hsu2012, Pan2025} and their references.  

This approach has two main limitations. First, that the growth and decline terms $f_g$ and $f_d$ must be restricted to a small set of functions for which the corresponding dynamical system may be solved explicitly; and second, that the transition between seasons must be sharp, disregarding in this way smoother types of transitions like those described by differentiable periodic reaction terms, commonly appearing in epidemiological models.

In this paper, we propose a methodology that allows the consideration of general growth and decline functions satisfying certain qualitative assumptions, as well as the treatment of  both the sharp and smooth transition cases. 
Our results are based on the \emph{bifurcation from an eigenvalue} technique introduced by Crandall and Rabinowitz in \cite{Crandall1971, Crandall1973}, and later extended in a series of papers. In our approach, this technique is initially applied to the smooth transition case and then extended to the sharp transition case by a limit process on a suitable approximating parameter.

In a nutshell, the idea of this technique is the following. Let $F(\tau,u)=0$ represent the equation \fer{prob}. An equilibrium of this problem is a function $u$ satisfying the equation in a single season, say $t\in[0,1]$, and such that $u(0)=u(1)$, so that when taken as the initial datum of the original problem, the solution is $1$-periodic. Suppose that $u\equiv0$ is an equilibrium for all $\tau$  and consider the linearization of problem \fer{prob} around zero, which may be expressed in terms of the Jacobian matrix $L(\tau):=D_uF(\tau,0)$.

In our problem, the common situation is that $u = 0$ is a stable equilibrium for small values of $\tau$, implying that all the eigenvalues of $L(\tau)$ are negative. The Crandall–Rabinowitz theorem states that if, for some $\tau_*$, the linearized problem loses stability, i.e., $\det(L(\tau_*)) = 0$, and this loss of stability is due to a single eigenvalue of $L(\tau_*)$ becoming zero,  that is, $\dim(\ker L(\tau_*)) = 1$, then, under additional transversality conditions on the zero-crossing of this eigenvalue, a bifurcation occurs.

The Crandall–Rabinowitz theorem, being a consequence of the implicit function theorem, provides local information about the bifurcation structure near $\tau_*$, in particular, the local uniqueness of the bifurcating branch. However, it does not preclude the possibility of further bifurcations occurring at larger values of $\tau$. Since such bifurcations, if they arise, do so from non-constant and explicitly unknown equilibrium states, the conditions under which they occur are generally analytically intractable. Consequently, one must typically rely on numerical computations to demonstrate the existence of these secondary bifurcations under parameter regimes that can be well approximated numerically.

Let us finally mention that Crandall–Rabinowitz technique has been extended to more complex settings in later works, providing a detailed local information on the bifurcating branch under additional regularity requirements, see, for example, \cite{Liu2007,Liu2013,Liu2013b}.

The rest of the paper is organized a follows. We finish this introduction section with an heuristic description of a simple scalar example of problem  \fer{prob}. Section~\ref{sec.setting} is devoted to the general formulation of the problem and some of its elementary properties. In Section~\ref{sec.results} we present our main results, which are then applied to the Lotka-Volterra-Malthus example in Section~\ref{sec.LVM}, where we also perform some numerical experiments illustrating our analytical results. Finally, in Section~\ref{sec.conclusion} we give our conclusions.

\subsection{A simple example}\label{subsec.ML}
Let us look more closely to the model by considering a simple scalar example,  see \cite{Liu2007} for a deeper analysis. 

Set the growth function in the logistic form $f_g(s) = r s(1-s)$ and the declining function as the Malthusian function $f_d(s) = -\mu s$, where $r,\mu>0$. For the extreme cases of only growth ($\tau=1$) and only decline ($\tau=0$), it is clear that the corresponding global stable equilibria are given by $u_g^* =1$ and $u^*_d =0$, respectively. 
Also notice that  $u^*_d =0$ is an equilibrium for all $\tau$, due to the property  $f_g(0)=f_d(0)=0$.
Then, the question of what equilibrium, if it does exist, may be the one corresponding to any $\tau\in(0,1)$ arises. 

\begin{figure}[t]
\centering
\includegraphics[width=0.32\linewidth]{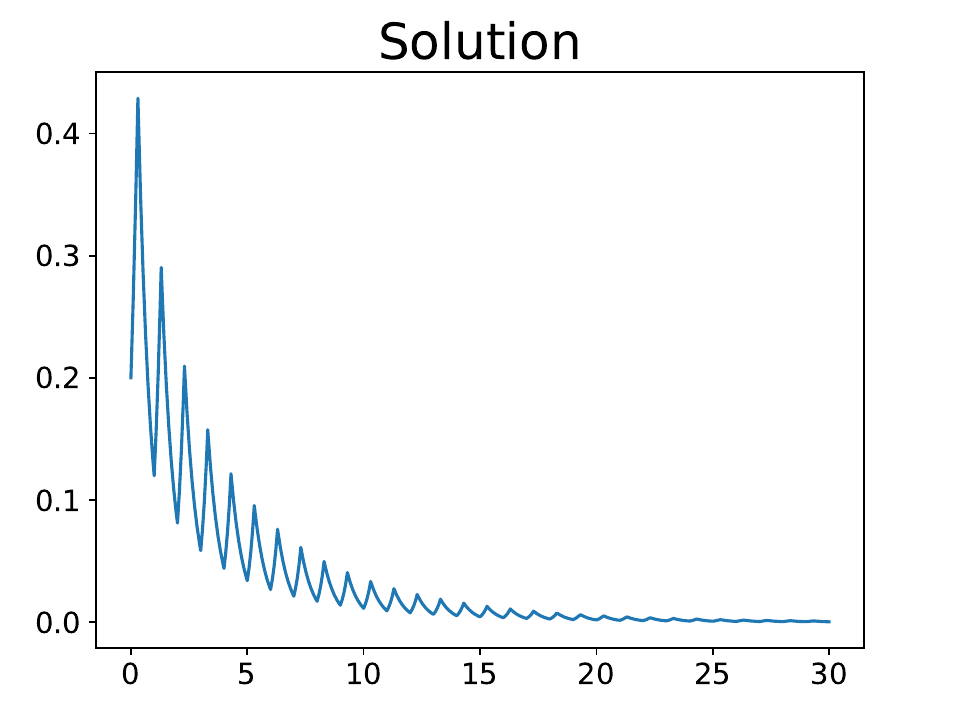}
\includegraphics[width=0.32\linewidth]{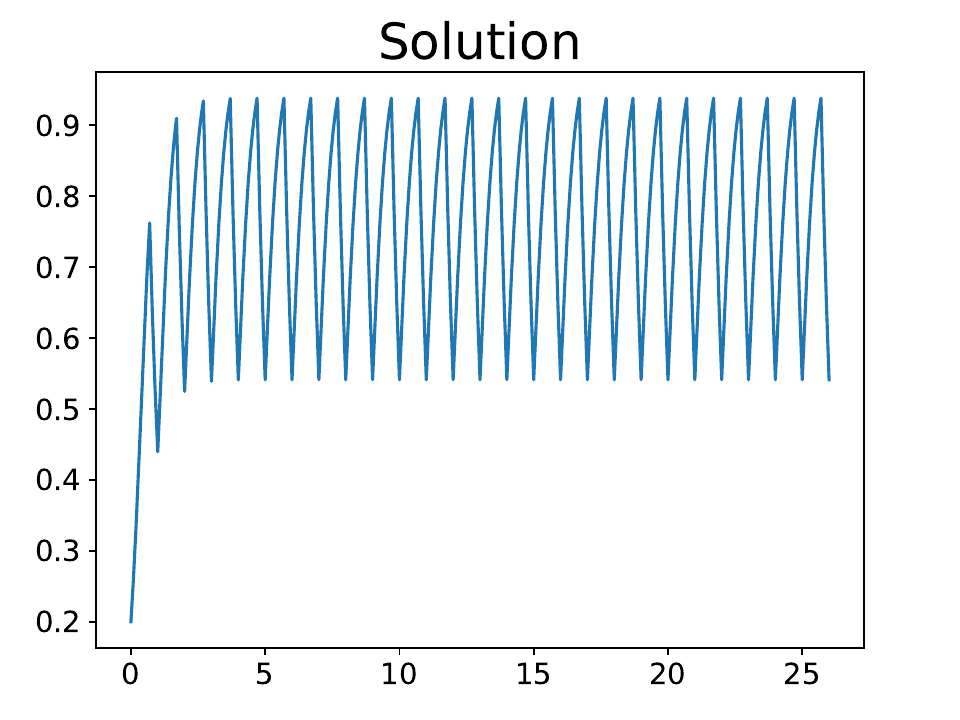}
\includegraphics[width=0.32\linewidth]{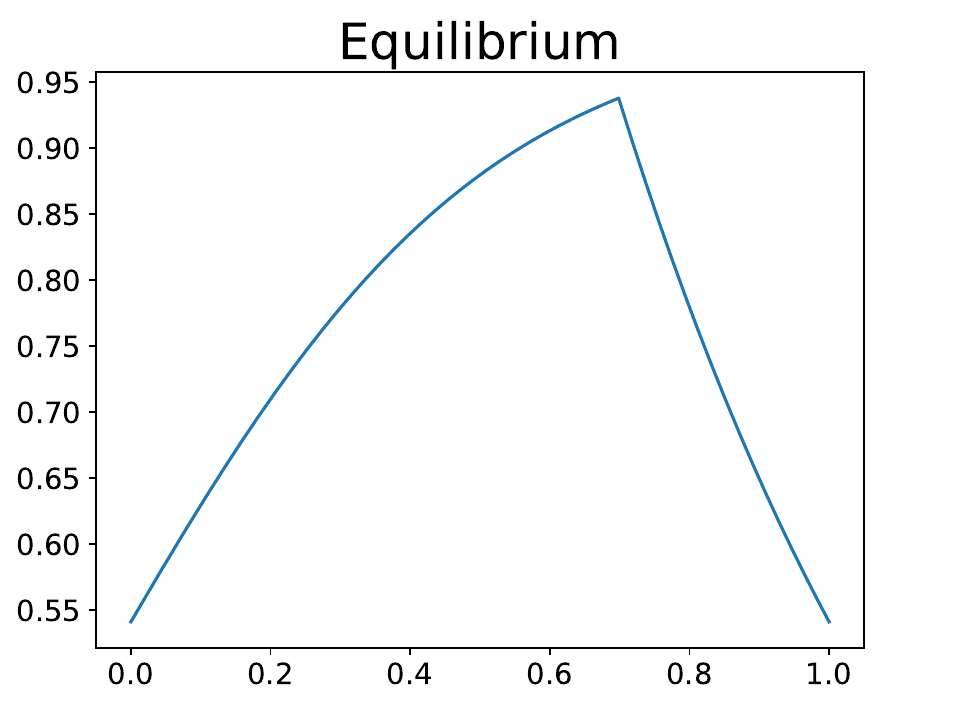}
  \caption{{\small Solution of problem \fer{prob} with Logistic growth and Malthusian decline. Left: Solution for $\tau=0.3$. Center: Solution for $\tau=0.7$. Right: Equilibrium for $\tau=0.7$. Other parameters are fixed as $u_0=0.2$, $r=0.01$ and $\mu=0.005$.}} 
\label{fig.intro1}
\end{figure}

Intuitively, it seems reasonable that if the growth period is short enough, this is, $\tau$ is close to zero, the extinction equilibrium will be the equilibrium of the new situation since the population density growth achieved during the growth season will not be enough to compensate the population decrease produced during the decline season, see Fig.~\ref{fig.intro1}, left. 

Reasoning in the same manner, it could be predicted that the global equilibrium for $\tau$ close to one  -mostly growth- should be that of the environment carrying capacity $u^*=1$ or by some quantity close to it. However,  this is prevented by the property $f_d(s)\neq 0$ for any $s\neq0$ so that, if a new equilibrium does exist it cannot be given by a constant steady state, see Fig.~\ref{fig.intro1}, center.  

The non-zero equilibrium of this dynamics must be a $1-$periodic function such that $u(n)=u(n+1)$, for $n=0,1,\ldots$, so that the population increase during  the growth season is exactly compensated by the population decrease during the declining season, see Fig.~\ref{fig.intro1}, right. Since the duration of these seasons is ruled by $\tau$, the  obvious conjecture is that $\tau$ should act as a bifurcation parameter separating two behaviors of the solution: declining to extinction, if $\tau<\tau_*$ and tending to the $1-$periodic non-zero equilibrium if $\tau>\tau_*$, where $\tau_*\in(0,1)$ is the critical bifurcation parameter.

\section{General setting of the problem and definitions}\label{sec.setting}

The model \fer{prob} considered in the Introduction captures a sharp, discontinuous transition between seasons which is not always appropriate. Indeed, in many natural phenomena such as  climatic seasons, transitions are smooth and progressive in time.  For capturing this effect we consider a regularized model in which this smoothing is taken into account and where the sharp transition model is recovered as a limit process, see Remark~\ref{convergence}.

Let $A\subset\R$ be a measurable set and denote by $\chil_A^0(t)$ to its characteristic function which is, in consequence, a measurable function.
 Let $\rho\in C^\infty(\R)$  be supported in $(-1/2,1/2)$ and satisfy
\begin{align}\label{hip.kernel}
\rho(s) \geq 0,\quad \rho'(s)s\geq 0,\quad \rho(s)=\rho(-s),\quad \int_\R \rho(s)ds = 1.	
\end{align}
Then, for $\eps\in(0,1)$, the rescaled kernel   $\rho_\eps(s) = \rho(s/\eps)/\eps$ is a Dirac $\delta-$approximation with support contained in $(-\eps/2,\eps/2)\subset (-1,1)$. We define 
\begin{align*}
\chil_{A}^\eps =	\rho_\eps * \chil_{A}.
\end{align*}
  Notice that, by standard approximation results, see e.g. \cite[Th. 4.22]{Brezisbook},  if $A$ is of finite measure then 
  $\chil_{A}^\eps \to \chil_A^0$ as $\eps\to0$ in $L^p(\R)$, for any $1\leq p<\infty$.

We are now prepared to give the general formulation of our model.
Let the season length be given by  $\tau\in(0,1)$. We introduce the partition of the time domain  $	\R_+ = G \cup D$ with $G=\bigcup_{n=0}^\infty G_n$, $D=\bigcup_{n=0}^\infty D_n$, and with $G_n,D_n$ given by \fer{def.cndn}. 

Consider two  functions 
\begin{align}
\label{ass.lip}
	f_g,f_d \in C^1(\R^N;\R^N)
\end{align}
 and an initial data $u_0\in \R^N$. Let $\eps\in[0,1)$. The problem is: Find $u:[0,\infty)\to\R^N$ such that $u(0)=u_0$, and
\begin{align}
&u'(t) = 	f_g(u(t)) \chil_G^\eps(t) + f_d(u(t)) \chil_D^\eps(t) && \text{for }t>0,&\label{prob2}\\
& u(0)=u_0. &&&\label{prob2b}
\end{align}
Observe that this formulation encapsulates both the sharp transition model ($\eps=0$) and the smooth transition model ($\eps>0$).

We assume
\begin{align}
&f_g(0)=f_d(0)=0,\label{ass.cero1}
\end{align}
implying that $u_*=0$ is an equilibrium for all $\tau\in(0,1)$, and the existence of a constant $K\in\R^N$ with positive components such that for $j=1,\ldots,N$ and $ u_j\in[0,K_j] $ we have 
\begin{align}
	f_{g,i}(u)\geq0,~f_{d,i}(u)\geq0 \text{ if } u_i=0,\quad  	f_{g,i}(u)\leq0,~f_{d,i}(u)\leq0 \text{ if } u_i=K_i,\label{ass.cero2}
\end{align}
implying the uniform boundedness of the solution of \fer{prob2}-\fer{prob2b}, this is,  if for $i=1,\ldots,N$ we have $u_{0,i}\in[0,K_i]$  then $u_{i}(t)\in[0,K_i]$ for $t\in[0,\infty)$.

For $\eps>0$, the existence and uniqueness of a $C^{2}((0,\infty))\cap C([0,\infty))$ solution and  its smooth dependence on the parameters  is ensured by classical results for  non-autonomous systems of differential equations.

For $\eps=0$, the solution to \fer{prob2}-\fer{prob2b} may be constructed by matching the solutions of the collection of successive  regular problems defined in $G_n$ and $D_n$, and then the existence, uniqueness and Lipschitz continuity of the solution is a straightforward consequence of the differentiability of $f_g$ and $f_d$.  The non-negativity of the components and the equilibria and boundedness properties also follow from \fer{ass.cero1}-\fer{ass.cero2}.

We conclude this section by introducing the notion of equilibrium we use. 
Let us introduce the notation $\chil_g^\eps(t,\tau) =\chil^\eps_{G_0}(t)$,  and  $\chil_d^\eps(t,\tau) =\chil^\eps_{D_0}(t)$, this is,
\begin{align}\label{def.chis2}
& \chil_g^\eps(t,\tau) = \int_0^\tau 	 \rho_\eps(t-s)ds, \quad
 \chil_d^\eps(t,\tau) =\int_\tau^1 	 \rho_\eps(t-s)ds,
\end{align}
and consider the increasing function 
\begin{align}\label{def.calR}
r_\eps(\tau) = \int_0^\tau\int_0^1 \rho_\eps(t-s)dt ds,
\end{align}
for which we have $r_\eps(\tau)\to\tau$ as $\eps\to0$, uniformly in $\tau$.

\begin{definition}\label{def.eq}
Let $\eps,\tau\in(0,1)$. The function $u^*\in   C^2_{per}([0,1])$ is an \emph{equilibrium} of \fer{prob2}-\fer{prob2b} if it satisfies $\abs{u^*_i(t)}\leq K_i$ for $i=1,\ldots,N$, and 
\begin{align*}
&(u^*)'(t)  = 	f_g(u^*(t)) \chil_g^\eps(t,\tau) + f_d(u^*(t)) \chil_d^\eps(t,\tau)\big),&&\text{for } t\in(0,1). 
\end{align*}
For $\eps=0$ the definition is similar, but with $u^*\in  W^{1,\infty}((0,1))\cap C_{per}([0,1])$.

\end{definition}

\begin{remark}\label{convergence}

The boundedness property assumed on the definition of an equilibrium is motivated by the  boundedness of the solution of \fer{prob2}-\fer{prob2b} implied by assumption \fer{ass.cero2}. Notice that this property together with Arzel\`a-Ascoli theorem ensure  the uniform convergence of any sequence of equilibria $u^*_\eps$ to an equilibrium $u^*_0$ as $\eps\to0$, where $u^*_\eps$ and $u^*_0$ denote equilibria of the problems corresponding to $\eps>0$ and $\eps=0$, respectively. 	
\end{remark}

\section{Main results}\label{sec.results}

Consider the Banach spaces  $X= C^1((0,1))\cap C_{per}([0,1]) $, $Y= C([0,1])$. The application of the Crandall-Rabinowitz theorem to our problem \fer{prob2}-\fer{prob2b} encompasses the analysis of the operator $F:(0,1)\times X \to Y$ given by
\begin{align}\label{def.F}
	F(\tau,v) = v' - \big(f_g(v) \chil_g^\eps(\cdot,\tau) + f_d(v) \chil_d^\eps(\cdot,\tau)\big). 
\end{align}
We also introduce some matrices that will play an important role in our analysis:
\begin{align*}
 & H[\tau,u](t) =	Df_g(u)\chil_g^\eps(t,\tau)+ Df_d(u) \chil_d^\eps(t,\tau),\\
& A[u]= Df_g(u)- Df_d(u). 
\end{align*}

Our main result establishes the conditions under which a uniquely determined local in time non-trivial branch of equilibria, $u^*$, emerges from the trivial equilibrium $u_*=0$.

\begin{theorem}\label{th.CR}
Assume \fer{hip.kernel}, \fer{ass.lip}, and \fer{ass.cero1}-\fer{ass.cero2}.
Let $u_*(t)=0$ for $t\in[0,1]$ and let $E\subset X$ be an open neighbourhood of $u_*$.	
Suppose that the following conditions are satisfied:
\begin{enumerate}
\item[(C1)]  
 There exists $\tau_*\in(0,1)$ such that the linear subspace of solutions of 
\begin{align}\label{prob.sys1}
\Phi'(t) =H[\tau_*,u_*](t)\Phi(t)\qtext{for } t\in(0,1),\quad \Phi(0)=\Phi(1)	,
\end{align}
is span by some non-zero function $\Phi_0$.
\item[(C2)] It is satisfied the transversality condition 
\begin{align}\label{cond.trans2}
\int_0^1 \rho_\eps(s-\tau_*) \Phi_R(s) \cdot A[u_*] \Phi_0(s) ds \neq 0,
\end{align}
where $\Phi_R$ spans the linear subspace of solutions of the dual problem
\begin{align}\label{prob.sys1d}
\Phi'(t) = -H[\tau_*,u_*](t)^T\Phi(t) \qtext{for } t\in(0,1), \quad \Phi(0)=\Phi(1)	.
\end{align}
\end{enumerate}
Then, there exist 
\begin{enumerate}
\item An open neighborhood $B\subset X$ of $u_*$ and an  interval $I=(-\delta,\delta)$, for some $\delta>0$,	
\item Functions $\tau^* \in C^1(I;(0,1))$, $u^*\in C^1(I;E)$,  with $\tau^*[0]=\tau_*$, $u^*[0]=u_*$ and $u^*[s]\neq u_*$ for $s\neq0$,
\end{enumerate}
 such that 
\begin{align*}
F(\tau^*[s],u^*[s]) =0 \qtext{for all }s\in I.	
\end{align*}
Moreover, for $N=\{(\tau,u)\in ( (0,1) \times E\backslash\{0\}) ) : F(\tau,u)=0 \}$, we have 
\begin{align*}
B \cap N = \{ (\tau^*[s], u^*[s]): s\in I\backslash \{0\}\}	.
\end{align*}
\end{theorem}

\begin{remark}
The local behavior of the non-trivial equilibrium  $u^*[s]$ may be approximated  as  $u^*[s]=u_*+s\Phi_0 + s \xi[s]$, where $\xi$ is a continuous function with $\xi[0]=0$.
Under further regularity assumptions on $f_g$ and $f_d$, a more detailed local behavior of $\tau^*[s]$ may  also be obtained, see \cite[Th.~1.18]{Crandall1973}.
\end{remark}

Conditions (C1) and (C2) are easy to check when the Jacobian matrices $Df_g(u_*)$ and $Df_d(u_*)$ are symmetric and commute. In this case, $H[\tau_*,u_*](t)$ and $H[\tau_*,u_*](s)$ are also symmetric and commute for all $t,s\in(0,1)$, implying that the solutions of \fer{prob.sys1} and \fer{prob.sys1d} may be expressed in terms of 
\begin{align}\label{def.calH}
\cH[\tau_*,u_*](t) = \int_0^t H[\tau_*,u_*](s)ds	
\end{align}
 as $\Phi_0(t) = e^{\cH[\tau_*,u_*](t)}\vfi_0$, and $\Phi_R(t) = e^{-\cH[\tau_*,u_*](t)}\vfi_0$, with $\vfi_0\in\ker(I-e^{\cH[\tau_*,u_*](1)})$.

\begin{corollary}\label{cor.symm}
In the conditions of Theorem~\ref{th.CR}, suppose that  $Df_g(u_*)$ and $Df_d(u_*)$  are symmetric and commute, that  $A[u_*]$ is definite, and that there exists $\tau_*\in(0,1)$ such that (C1) is satisfied. Then (C2) is also satisfied and the thesis of the theorem applies.	
\end{corollary}

If $H[\tau_*,u_*](t)$ does not commute with itself for different times,  we must resort to the notion of fundamental matrix for solving multi-dimensional non-autonomous non-homogeneous linear problems 
\begin{align*}
& \Psi(t)' = H[\tau_*,u_*](t) \Psi(t) + \psi(t) && \text{for }t>0,& \\
& \Psi(0)=\Psi_0. &&& 
\end{align*}
In this case, the solution is given by 
\begin{align}\label{solnh}
\Psi(t) = G(t,0)\Psi_0 + \int_0^t G(t,s)\psi(s)ds,	
\end{align}
where the 
 non-singular matrix $G(t,s)$ is the unique solution of
\begin{align*}
\p_t G(t,t_0) = H[\tau_*,u_*](t) G(t,t_0)	
\end{align*}
such that $G(t_0,t_0)=I$, the identity matrix.
On noting that the fundamental matrix of the dual problem is $G(s,t)^T$, 
we deduce that the solutions of \fer{prob.sys1} and \fer{prob.sys1d} may be expressed as
\begin{align*}
\Phi_0(t) = G(t,0)\vfi_0,\quad 	\Phi_R(t) = G(1,t)^T\vfi_R, \qtext{for }t\in(0,1),
\end{align*}
with $\vfi_0\in\ker(I-G(1,0))$ and $\vfi_R\in\ker(I-G(1,0)^T)$. Then, the transversality condition \fer{cond.trans2} may be rewritten as 
\begin{align*}
\int_0^1 \rho_\eps(s-\tau_*) ~ \vfi_R^T G(1,s)AG(s,0)\vfi_0 ds \neq 0.
\end{align*}
It is not difficult to find examples in which this transversality condition is not satisfied, showing that it is a restrictive condition.

\bigskip

\no\emph{Proof of Theorem~\ref{th.CR}.} The proof consists on checking that $F$ defined by \fer{def.F} satisfies the four conditions of Crandall-Rabinowitz theorem \cite[Theorem~1.7]{Crandall1971}. The first two conditions are straightforward.

\emph{Condition (a): $F(\tau,u_*)=0$ for all $\tau\in(0,1)$. } This is clear, since $f_g(u_*)=f_d(u_*)=0$ by hypothesis \fer{ass.cero1}.

\emph{Condition (b): The partial Frechet derivatives $F_\tau,~F_v$ and $F_{\tau v}$ exist and are continuous. }Replacing \fer{def.chis2} in \fer{def.F} and differentiating with respect to $\tau$, we get, for $\sigma\in(0,1)$ 
\begin{align*}
F_\tau(\tau,v)[\sigma] = - (f_g(v)-f_d(v)) \rho_\eps(\cdot-\tau)\sigma .	
\end{align*}
We also have, for $\vfi\in X$,
\begin{align*}
F_v(\tau,v)[\vfi] & = \frac{d}{dh} F(\tau, v+h\vfi)\Big|_{h=0} \\
& = \vfi'-	\big(Df_g(v) \chil_g^\eps(\cdot,\tau) + Df_d(v) \chil_d^\eps(\cdot,\tau)\big) \vfi  ,
\end{align*}
and then 
\begin{align*}
F_{\tau v}(\tau,v)[\sigma,\vfi] = - (Df_g(v)-Df_d(v)) \rho_\eps(\cdot-\tau)\sigma \vfi.
\end{align*}
The continuity of these operators follows from the continuous differentiability of $f_g,~f_d$ and $\rho_\eps$, see \fer{hip.kernel} and \fer{ass.lip}.

\emph{Condition (c): $\ker(F_v(\tau_*,u_*))$ and $Y\backslash R(F_v(\tau_*,u_*))$  are one-dimensional. }
We have $\Phi\in \ker(F_v(\tau_*,u_*))$ if $\Phi$ is a solution of \fer{prob.sys1}. 
By hypothesis, the linear space of solutions of \fer{prob.sys1} is one-dimensional, i.e, the dimension of  $\ker(F_v(\tau_*,u_*))$ is one.

Notice that using formula \fer{solnh} with $\psi=0$ we obtain the representation of the solution of  \fer{prob.sys1} as $\Phi_0(t)=G(t,0)\vfi_0$, whenever the periodic boundary condition $(I-G(1,0))\vfi_0=0$ is satisfied. Since, by hypothesis, the solution is unique up to a multiplicative constant, we deduce $\dim(\ker(I-G(1,0)))=1$.

Now we check that the co-dimension of $R(F_v(\tau_*,u_*))$ is one. 
We have $\psi \in  R(F_v(\tau_*,u_*))$ if the equation 
\begin{align}
&	 \xi'(t) =	H[\tau_*,u_*](t) \xi(t) + \psi(t) &&\text{ for }t\in(0,1), \label{nR.1}
\end{align}
has a solution $\xi \in X$, this is, also satisfying $\xi(0)=\xi(1)$.  Let $\Phi_R$ be  the solution of the dual problem \fer{prob.sys1d}. Then, multiplying \fer{nR.1} by $\Phi_R$ and integrating in $(0,1)$ we obtain
\begin{align*}
-\int_0^1 \xi(t)\cdot \Phi_R'(t)dt =	\int_0^1 H[\tau_*,u_*](t)\xi(t)\cdot \Phi_R(t) + \int_0^1 \psi(t)\cdot \Phi_R(t) dt,
\end{align*}
 and using that $\Phi_R'(t) = -H[\tau_*,u_*](t)^T\Phi_R(t)$, we deduce the orthogonality condition
\begin{align*}
\int_0^1  \psi(t)  \cdot \Phi_R(t)dt = 0	.
\end{align*}
Therefore, any $\psi\in R(F_v(\tau_*,u_*))$ must be orthogonal  to the subspace span by $ \Phi_R(t)$, implying that the co-dimension of $R(F_v(\tau_*,u_*))$ is one.

\emph{Condition (d): $F_{\tau v}(\tau_*,u_*)[1,\Phi_0] \notin R(F_v(\tau_*,u_*))$. } This condition reads
\begin{align*}
0 \neq \int_0^1 \Phi_R(s) \cdot F_{\tau v}(\tau_*,u_*)[1,\Phi_0](s)  ds 
 =
- 	\int_0^1 \rho_\eps(s-\tau_*) \Phi_R(s) \cdot A[u_*] \Phi_0(s) ds,
\end{align*}
which is true by assumption.
\hfill$\Box$


\bigskip

\no\emph{Proof of Corollary~\ref{cor.symm}. }If the Jacobian matrices $Df_g(u_*)$ and $Df_d(u_*)$ are symmetric and commute then the matrices $A[u_*]$, $H[\tau_*,u_*](t)$ and $H[\tau_*,u_*](s)$ have the same property for all $t,s\in(0,1)$. Thus,  the solutions of \fer{prob.sys1} and \fer{prob.sys1d} may be expressed in terms of the matrix $\cH[\tau_*,u_*](t)$ given by \fer{def.calH} as 
\begin{align*}
\Phi_0(t) = e^{\cH[\tau_*,u_*](t)}\vfi_0,\quad \Phi_R(t) = e^{-\cH[\tau_*,u_*](t)}\vfi_R, 	
\end{align*}
with $\vfi_0\in \ker(I- e^{\cH[\tau_*,u_*](1)})$, one-dimensional by assumption,  and $\vfi_R\in \ker(I- e^{-\cH[\tau_*,u_*](1)})$. Since $e^{\cH[\tau_*,u_*](1)}$ is non-singular, we have $ \ker(I- e^{\cH[\tau_*,u_*](1)})= \ker(I- e^{-\cH[\tau_*,u_*](1)})$, so that we may choose $\vfi_R=\vfi_0$. On noting that  $A[u_*]$ also  commutes with $e^{\cH[\tau_*,u_*](t)}$, we deduce 
\begin{align*}
	 \Phi_R(t) \cdot A[u_*] \Phi_0(t) = \vfi_R^Te^{-\cH[\tau_*,u_*](t)}A[u_*]e^{\cH[\tau_*,u_*](t)}\vfi_0 = \vfi_0^T A[u_*]\vfi_0,
\end{align*}
 for all $t\in(0,1)$.  Finally, since by assumption $A[u_*]$ is definite, we have  $ \vfi_0^T A[u_*]\vfi_0\neq0$. Then, condition \fer{cond.trans2} in (C2) reduces to 
\begin{align*}
\int_0^1 \rho_\eps(s-\tau_*) ds \neq 0,
\end{align*}
which is true in view of $\tau_*\in(0,1)$ and the properties assumed on the kernel $\rho_\eps$, see \fer{hip.kernel}. 
\hfill$\Box$

\subsection{Secondary bifurcations}\label{rem.secbif}

It is common in systems of equations for the first bifurcation to emerge from the trivial equilibrium $u_*$, leading to a new equilibrium where only one component is non-zero. As $\tau$ increases further, additional bifurcations may occur in which other components of the equilibrium also become non-zero. We refer to these subsequent bifurcations as \emph{secondary bifurcations}. 

Crandall-Rabinowitz's theory may be still applied to this situation, although the modification of the operator $F$ that must be introduced to express its dependence on the emerging branch of non-trivial equilibria $u^*[s]$ leads to conditions that, in general, can not be checked analitically.

We give here a brief account of the conditions ensuring the emergence of secondary bifurcations.

Suppose that a primary bifurcation has emerged under the situation of Theorem~\ref{th.CR}, and consider the modification of the operator $F$ defined in \fer{def.F} given by 
\begin{align*}
	F(s,v) = v' +  f(\cdot,\tau^*[s], u^*[s]) -  f(\cdot,\tau^*[s], v+u^*[s])  ,
\end{align*}
for  $s\in I$, where 
 $f:[0,1]\times(0,1)\times \R^N\to\R^N$ is given by 
\begin{align*}
	 f(t,\tau ,z)= f_g(z) \chil_g^\eps(t,\tau) + f_d(z) \chil_d^\eps(t ,\tau).
\end{align*}

Set $v_*=0$. Condition (a) is trivially satisfied: $F(s,v_*)=0$ for all $s\in I$.  Condition (b) is also straightforward to check if we assume $f_g, f_d \in C^2(\R^N)$. In particular, we have 
\begin{align*}
F_v(s,v)[\vfi] & = \vfi'-	Df(\cdot,\tau^*[s], v+u^*[s]) \vfi  ,\quad \vfi\in X,
\end{align*}
where $Df$ is the matrix  $Df = \big(\dfrac{\p f_i}{\p z_j}\big)$ for $i,j=1,\ldots,N$, 
and 
\begin{align*}
F_{s v}(s,v)[\sigma,\vfi] = - B[s,v]\vfi \sigma	,\quad \sigma\in  \R, ~\vfi\in X,
\end{align*}
 where
\begin{align*}
B_{ij}[s,v] =& \frac{\p}{\p s}\frac{\p f_i}{\p z_j}(\cdot,\tau^*[s], v+u^*[s]) 
 = \frac{\p^2f_i}{\p \tau \p z_j}(\cdot,\tau^*[s], v+u^*[s])  \frac{\p \tau^*}{\p s}[s] \\
 &+
\sum_k \frac{\p^2f_i}{\p z_k \p z_j}(\cdot,\tau^*[s], v+u^*[s]) \frac{\p u_k^*}{\p s}[s].
\end{align*}

Condition (c): 
Let $G[s](t,\hat t)$ be the fundamental matrix corresponding to 
\begin{align}\label{eq.sec22}
\vfi'(t) =	Df(t,\tau^*[s], u^*[s](t)) \vfi(t),\quad t\in[0,1].
\end{align}
 We have  $\ker(F_v(s_*,v_*))$  one-dimensional if
there exists  $s_*\neq0$ such that  $\dim(\ker(I-G[s_*](1,0)))=1$. 
If such $s_*$ does exist and $\vfi_0$ spans this kernel, then $\Phi_0[s_*](t) = G[s_*](t,0)\vfi_0$ spans  $\ker(F_v(s_*,v_*))$ and it is thus one-dimensional.

For showing that the co-dimension of $R(F_v(s_*,v_*))$ is one we again resort to the dual problem. Like in the proof of Theorem~\ref{th.CR} we deduce that the space span by any solution $\Phi_R[s_*]$ of the dual problem is orthogonal to  $R(F_v(s_*,v_*))$, implying its one-dimensional property.

Condition (d). Again, like in the proof of Theorem~\ref{th.CR}, this condition takes the form
\begin{align*}
 \int_0^1 \Phi_R[s_*](t) \cdot F_{s v}(s_*,v_*)[1,\Phi_0[s_*]](t)  dt \neq 0,
 \end{align*}
or, in terms of the fundamental matrix
\begin{align}\label{int.cond22}
 \int_0^1 \vfi_0^T G[s_*](1,t) B[s_*,v_*] G[s_*](t,0)\vfi_0  dt \neq 0.
 \end{align}
 
 Notice that the verification of conditions (c) and (d) depend on the knowledge of the analytical form of the primary bifurcation branch $u^*[s]$ for $s\in I$, which is not available in general. 

\section{Competitive Lotka-Volterra growth versus Malthusian decline}\label{sec.LVM}

In this section, we analyze an example of a seasonal dynamical system in which both primary and secondary critical bifurcation parameters can be computed analytically. For the primary bifurcation, the Jacobian matrices $Df_g(u_*)$ and $Df_d(u_*)$ are constant and diagonal, allowing the application of Corollary~\ref{cor.symm} to determine the primary critical bifurcation parameter $\tau_*$. 

For assesing the possible emergence of a secondary bifurcation, we study Condition (c) of the theorem, this is, the existence of some $s_*$ such that  $\dim(\ker(I-G[s_*](1,0)))=1$, where $G$ is the fundamental matrix corresponding to \fer{eq.sec22}, and thus depends on the primary bifurcation branch of equilibria $u^*[s]$. 

Approximating such $s_*$ may be done by numerical methods that involve the computation of  the equilibria $u^*[s]$. In our particular example, $s_*$ may also be computed  analitically, allowing us to checking the validity of the numerical approximation.

To ensure the actual emergence of a secondary bifurcation Condition (d) should also be checked. However, the analytical form of the integral \fer{int.cond22} expressing this condition is not at hand.

Like in \cite{Koch1974, Namba1984, Pan2025}, we consider problem \fer{prob2}-\fer{prob2b} for two populations ruled by the competitive Lotka-Volterra model in the growth season,
\begin{align*}
f_{g,i}(s_1,s_2) = s_i(\alpha_i -\beta_{i1}s_1-\beta_{i2}s_2),
\end{align*}
and by the Malthusian model in the decline season,
\begin{align*}
f_{d,i}(s_1,s_2) = -\mu_i s_i	,
\end{align*}
for $i=1,2$, with $\mu_i>0$. We further assume the following conditions on the Lotka-Volterra coefficients: $\alpha_i,\beta_{ii}> 0$,  $\beta_{ij}\geq0$, and
\begin{align*}
\beta_{ii}\alpha_j >\beta_{ji}\alpha_i,\quad \beta_{11}\beta_{22}>\beta_{12}\beta_{21},	
\end{align*}
for $i,j=1,2$, $i\neq j$, ensuring that the coexistence equilibrium is stable for the only-growth dynamical system.

 \paragraph{Primary bifurcation} Let $u_*=0$ be the trivial equilibrium corresponding to this functional choice for problem \fer{prob2}-\fer{prob2b}. Since the Jacobian matrices
\begin{align*}
 Df_g(u_*) = \begin{pmatrix}
 \alpha_1 & 0 \\ 0 & \alpha_2	
 \end{pmatrix},\quad
 Df_d(u_*) = \begin{pmatrix}
 -\mu_1 & 0 \\ 0 & -\mu_2	
 \end{pmatrix},
\end{align*}
are diagonal and the matrix $A= Df_g(u_*)- Df_d(u_*)$ is diagonal and definite positive, we are in the conditions of Corollary~\ref{cor.symm}. 
Thus, we only have to check if Condition (C1) of Theorem~\ref{th.CR} is satisfied, this is, whether a $\tau_*\in(0,1)$ does exist such that $\ker(I-G(1,0))$ is one-dimensional.

The solution of the linearized equilibrium problem \fer{prob.sys1} admits a representation through the exponential matrix $G(t,0)=e^{\cH[\tau,u_*](t)}$,  with
\begin{align*}
\cH[\tau,u_*](1) 
& = \begin{pmatrix}
(\alpha_1+\mu_1) r_\eps(\tau) -\mu_1 r_\eps(1) & 0 \\ 0 &  (\alpha_2+\mu_2) r_\eps(\tau) -\mu_2 r_\eps(1)	
 \end{pmatrix}.	
\end{align*}
Since $r_\eps$ is increasing with $r_\eps(0)=0$, see \fer{def.calR}, each one of the equations
\begin{align}\label{eq.tau}
r_\eps(\tau) = 	\frac{\mu_i}{\alpha_i+\mu_i} r_\eps(1) \qtext{for }i=1,2,
\end{align}
has a solution $\tau_i\in(0,1)$, for $i=1,2$. Assume $\alpha_1/\mu_1\neq\alpha_2/\mu_2$ so that $\tau_1\neq\tau_2$ and suppose, without loss of generality, that $\alpha_1/\mu_1>\alpha_2/\mu_2$, implying $\tau_*:=\tau_1<\tau_2$. Then $\dim(\ker(I-e^{\cH[\tau_*,u_*](1)}))=1$, and condition (C1) of Theorem~\ref{th.CR} is satisfied.
Therefore, a bifurcation occurs at $\tau_*$. 

For computing the approximation  $ u^*[s](t) \approx u_*+s\Phi_0(t)$ to the new equilibrium, we notice that
\begin{align*}
I - e^{\cH[\tau_*,u_*](1)} = \begin{pmatrix}
0 & 0 \\ 0 &  1 - a	
 \end{pmatrix},	
\end{align*}
where $a=\exp((\alpha_2+\mu_2) r_\eps(\tau_*) -\mu_2 r_\eps(1))<1$. Therefore, 
$\ker(I-e^{\cH[\tau_*,u_*](1)})$ is span by $\vfi_0=(1,0)$, and the approximation is given by
\begin{align*}
& u^*[s](t) \approx u_* +  s\begin{pmatrix} e^{\cH_{11}[\tau_*,u_*](t)}\\ 0 \end{pmatrix},\\
& \text{with }
\cH_{11}[\tau_*,u_*](t) = \int_0^t (\alpha_1\chil_g^\eps(z,\tau_*) - \mu_1  \chil_d^\eps(z,\tau_*))dz	,
\end{align*}
for $s$ small.

 \paragraph{Secondary bifurcation} 
  We check the existence of some $s_*$ such that  $\dim(\ker(I-G[s_*](1,0)))=1$, where $G[s]$ is the fundamental matrix corresponding to \fer{eq.sec22}, with 
  \begin{align*}
  Df(t,\tau^*[s], u^*[s](t)) = 	Df_g(u^*[s](t)) \chil_g^\eps(t,\tau^*[s])+ Df_d(u^*[s](t)) \chil_d^\eps(t,\tau^*[s]),
  \end{align*}
being 
\begin{align*}
 &Df_g(u^*[s](t)) = \begin{pmatrix}
 \alpha_1-2\beta_{11} u^*_1[s](t) & -\beta_{12}u^*_1[s](t) \\ 0 & \alpha_2  -\beta_{21}u^*_1[s](t)	
 \end{pmatrix},\\
 &
 Df_d(u^*[s](t)) = \begin{pmatrix}
 -\mu_1 & 0 \\ 0 & -\mu_2	
 \end{pmatrix}.
\end{align*}
We consider separately the diagonal and the non-diagonal cases.

\emph{Diagonal case. } If $\beta_{12}=0$,  $G[s]$ may be expressed in exponential form so that  we have to check whether  $s_*\in I$ does exist such $\ker\big(I-\exp(\cH[s_*](1))\big)$ is one-dimensional, where $\cH[s](1)$ is a diagonal matrix with
\begin{align*}
&\cH_{11}[s](1) = (\alpha_1 +\mu_ 1)r_\eps(\tau^*[s])-2\beta_{11}U[s]- \mu_1 r_\eps(1), \\
&\cH_{22}[s](1) = (\alpha_2 +\mu_ 2)r_\eps(\tau^*[s]) - \beta_{21}U[s] - \mu_2 r_\eps(1), 
\end{align*}
being
\begin{align}\label{def.I}
U[s] = 	 \int_{0}^{\tau^*[s]}  u^*_1[s](t)dt>0 .
\end{align}
In other words, we must solve the equation for $s_*$ given by 
\begin{align}\label{eq.tau2}
(\alpha_2 +\mu_ 2)r_\eps(\tau^*[s_*]) - \beta_{21}U[s_*] = \mu_2 r_\eps(1).
\end{align} 
If a solution does exist, the secondary bifurcation is approximated by
\begin{align}
& u^{**}[\sigma](t) = u^*[s_*](t) + \sigma \begin{pmatrix} 0 \\ e^{ \cH_{22}[s_*](t)}	
 \end{pmatrix},\qtext{with } \label{sec.bif}\\
&\cH_{22}[s_*](t) = \int_0^t \big( (\alpha_2-\beta_{21}u_1[s_*](z))\chil_g^\eps(z,\tau^*[s_*]) - \mu_2  \chil_d^\eps(z,\tau^*[s_*])\big) dz	,\label{sec.bif2}
 \end{align}
for $\sigma$ small.

\emph{Non-diagonal case. }If $\beta_{12}\neq0$ the fundamental matrix $G[s]$ corresponding to \fer{eq.sec22} is given by $G_{21}[s](t,\hat t) =	0$,
\begin{align*}
G_{ii}[s](t,\hat t) =& 	\exp(\cH_{ii}[s](t)-\cH_{ii}[s](\hat t)),\qtext{for }i=1,2,\\
G_{12}[s](t,\hat t) =& 	 \exp(\cH_{11}[s](t)\\
& -\cH_{11}[s](\hat t))\int_{\hat t}^t  \exp(\cH_{22}[s](z)-\cH_{11}[s](z)) H_{12}[s](z) dz,
\end{align*}
where  $\cH_{22}[s]$ is given by \fer{sec.bif2}, and 
\begin{align*}
\cH_{11}[s](t) &= \int_0^t \big( (\alpha_1-2\beta_{11} u_1^*[s](z))\chil_g^\eps(z,\tau^*[s]) - \mu_1  \chil_d^\eps(z,\tau^*[s])\big) dz	,\\
\cH_{12}[s](t) &= -\int_0^t \beta_{12} u_1^*[s](z)\chil_g^\eps(z,\tau^*[s]) dz	.
\end{align*}
 Therefore, if there exists $s_*\in I$ solving equation \fer{eq.tau2}, then the direct and dual problems \fer{prob.sys1} and \fer{prob.sys1d} have solutions span by $\Phi_0(t)=G[s_*](t,0)\vfi_0$ and $\Phi_R(t)=G[s_*](1,t)^T\vfi_R$, where 
 \begin{align*}
 \vfi_0 =  \begin{pmatrix}
\dfrac{G_{12}[s_*](1,0)}{1-G_{11}[s_*](1,0)} 	\\[1em] 1
 \end{pmatrix}, \quad
  \vfi_R =  \begin{pmatrix}
0 	\\[1em] 1
 \end{pmatrix}.		
 \end{align*}
Notice that in this non-diagonal case the secondary bifurcation shares contributions of both components, as the following approximation indicates
\begin{align}
 u^{**}[\sigma](t) &= u^*[s_*](t) + \sigma \begin{pmatrix}
\vfi_{0,1} 	\exp(\cH_{11}[s_*](t)) + G_{12}[s_*](t,0)  \\ 	\exp(\cH_{22}[s_*](t))
 \end{pmatrix}, \label{sec.bifnonsym}
\end{align}
for $\sigma$ small.

\paragraph{Analytical computation}\label{rem.exac}

For our example (with $\eps=0$), after but close to the first bifurcation we have $u_1^*[s] >0$ while $u_2^*[s]\equiv 0$, for $s\in I$. Thus, the system of equations is uncoupled, and $u_1^*[s]$ solves the logistic-Malthusian equation
\begin{align*}
&	v'(t) = v(t) (\alpha_1-\beta_{11}v(t))& &\text{for }t\in (0,\tau^*[s]), & \\
&	v'(t) = -\mu_1 v(t) & &\text{for }t\in (\tau^*[s],1), & 
\end{align*}
for $s\in I$.
 Solving for periodic boundary conditions and integrating in $(0,\tau^*[s])$ yields
 \begin{align*}
 U[s] = \frac{1}{\beta_{11}}(\alpha_1 \tau^*[s] -\mu_1 (1-\tau^*[s])),
 \end{align*}
so that the solution of \fer{eq.tau2} is given by
\begin{align}\label{taustar2}
\tau^*[s_*] = \Big(1+\frac{\alpha_2\beta_{11}-\alpha_1\beta_{21}}{\mu_2\beta_{11}-\mu_1\beta_{21}}\Big)^{-1}. 	
\end{align}

\subsection{Numerical experiments}

For the numerical discretization of the dynamical system \fer{prob2}-\fer{prob2b}, we take an uniform mesh of size $\Delta t$ resulting in $\texttt{nper} =1+ 1/\Delta t$ points for each period (of length one) and $\texttt{T*nper}+1$ for the whole time interval $[0,T]$. The size $\Delta t$ is chosen to have $\texttt{nper}\in\N$.

Then we solve problem \fer{prob2}-\fer{prob2b} in succesive periods by an explicit Runge-Kutta method, as implemented in  \cite[solve\_ivp]{Scipy2020}. We consider that an equilibrium has been reached if the difference between two consecutive periods satisfies
\begin{align}\label{tol}
\norm{u_{n+1}-u_n}_{L^2(0,1)}< \texttt{tol},	
\end{align}
where \texttt{tol} is a tolerance parameter, see Table~\ref{table:1}.

\begin{table*} 
\small
{
\centering
\begin{tabular}{|c|c|c|c|c|c|c|}
\hline 
\hline
\multicolumn{7}{|c|}{$N=2,~\eps = 0$} \\ 
\hline
$T_{max}$& $\Delta t$ & \texttt{tol} & $u_0$ & $\alpha_1,\alpha_2$ & $\beta_{11},\beta_{21},\beta_{22}$ & $\mu_1,\mu_2$\\
\hline
50000   & 0.1/365  & 1.\,e-15  & $(1,0.25)$ & 2, 1  & 1, 0.25, 1 & 1, 1\\ \hline
\hline
\end{tabular}
\caption{\small Parameters that are fixed in all the experiments. $T_{max}$: maximum number of periods evaluated until reaching the equilibria. $\Delta t$: time step for the dynamical system discretization. \texttt{tol}: tolerance for equilibrium, see \fer{tol}. $u_0$: initial datum for the dynamical system. Its value is half the equilibrium of the only-growth Lotka-Volterra system. $\alpha_i,\beta_{ij}$: Lotka-Volterra coefficients. $\mu_i$: Malthusian coefficients.}
\label{table:1}
}
\end{table*}

In our experiments, there is not a relevant difference between the solutions corresponding to $\eps=0$ and $\eps>0$ small,  as expected in view of the convergence property commented in Remark~\ref{convergence}. Therefore, we only show the experiments performed with $\eps=0$, for which we have exact computations of the primary and secondary critical bifurcation parameters.

\subsubsection{Parameter choice}

For observing a first bifurcation from the trivial equilibrium we must set $\alpha_1/\mu_1 \neq \alpha_2/\mu_2$, which is the condition ensuring that the solutions of equations \fer{eq.tau} for $\tau$ are different, and thus the dimension of $\ker(I-G(1,0))$ is one. 
We thus set
\begin{align*}
	\alpha_1 = 2,\quad  \alpha_2 = 1,\quad \mu_1=\mu_2=1. 
\end{align*}
For the intra-competition coefficients we set $\beta_{11}=\beta_{22}=1$, and for the inter-competition coefficients we take $\beta_{21}=0.25$ and then distinguish the two cases deduced in the previous section for the secondary bifurcation: the diagonal case, $\beta_{12}=0$, and the non-diagonal case $\beta_{12}>0$.

With this parameter choice the primary bifurcation takes place at 
\begin{align*}
\tau_* = \frac{\mu_1}{\alpha_1+\mu_1} = \frac{1}{3}.	
\end{align*}

In our example, the secondary critical bifurcation parameter may also be computed exactly, see  \fer{taustar2}. For our parameter choice, we get $\tau^*=0.6$.

We also compute the value of $\tau^*$ provided by the numerical resolution of equation  \fer{eq.tau2}, that for $\eps=0$, takes the form 
\begin{align}\label{eq.tau2b}
(\alpha_2 +\mu_ 2)\tau^*[s_*] - \beta_{21} U[s_*] = \mu_2 .
\end{align}
For checking whether equation \fer{eq.tau2b} is satisfied in the numerical discretization,  we consider the mesh $\tau_n = n /365$, for $n = 0, \ldots, 365$, and  
 then verify if there exists an index $n^* \in \{0, \ldots, 365\}$ such that
\begin{align*}
(\alpha_2 + \mu_2)\tau_{n^*} - \beta_{21} \hat{U}[n^*] \approx \mu_2,
\end{align*}
where $\hat{U}[n]$ denotes a numerical approximation of the integral  \fer{def.I}. More precisely, it approximates the expression
\begin{align*}
 \int_{0}^{\tau_{n}}  u^*_1[n](t)dt,
\end{align*}
where $u^*[n]$ is the numerically computed equilibrium solution corresponding to $\tau = \tau_{n}$.

\begin{figure}[t]
\centering
\includegraphics[width=0.32\linewidth]{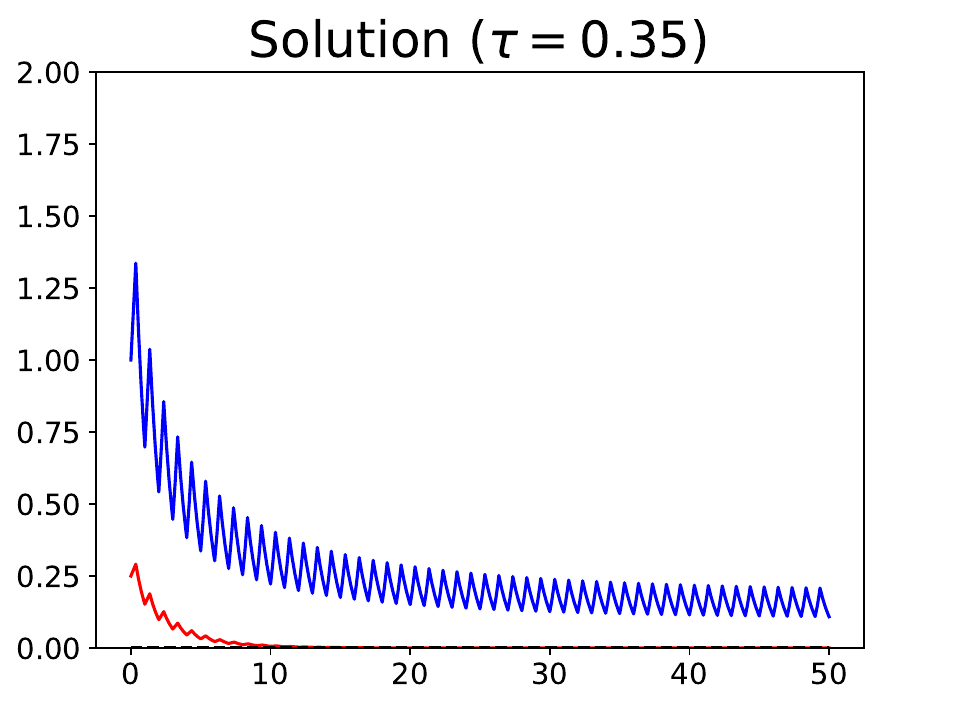}
\includegraphics[width=0.32\linewidth]{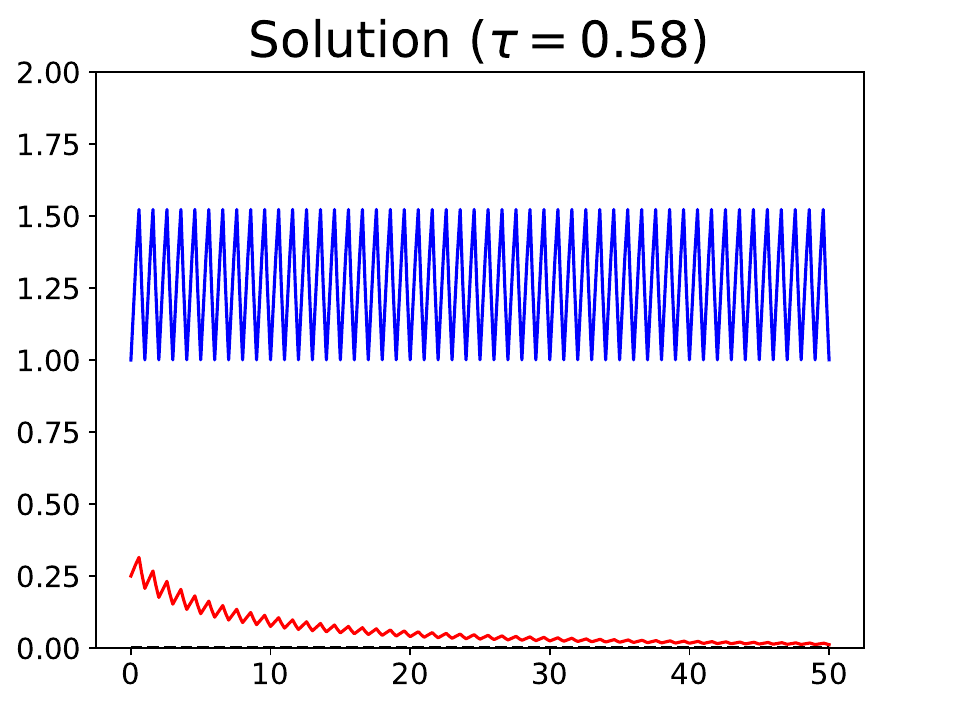}
\includegraphics[width=0.32\linewidth]{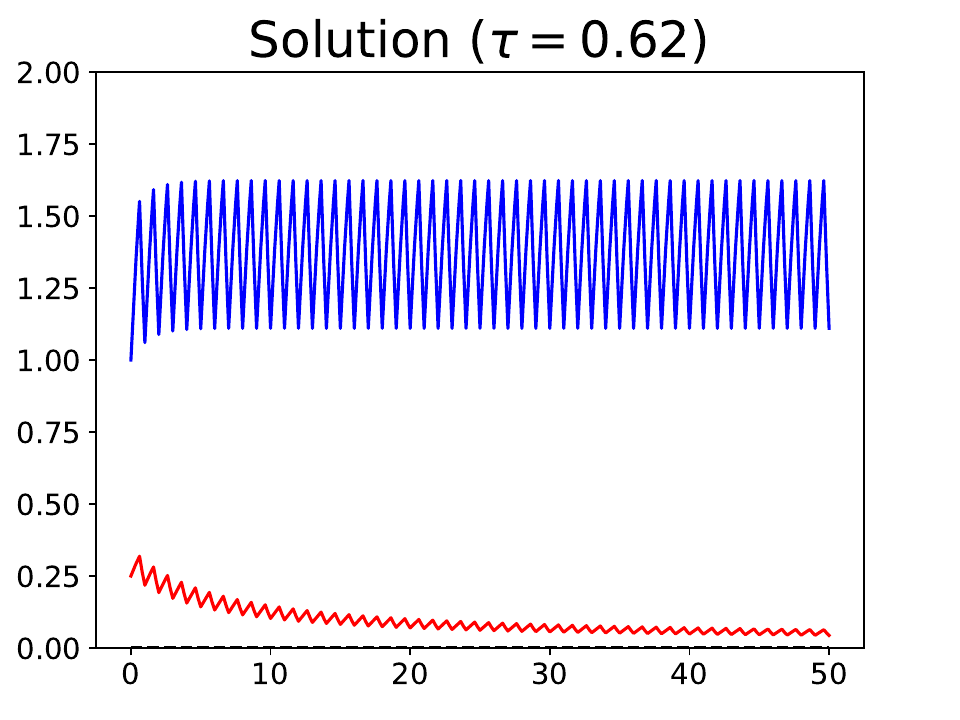}\\
\includegraphics[width=0.32\linewidth]{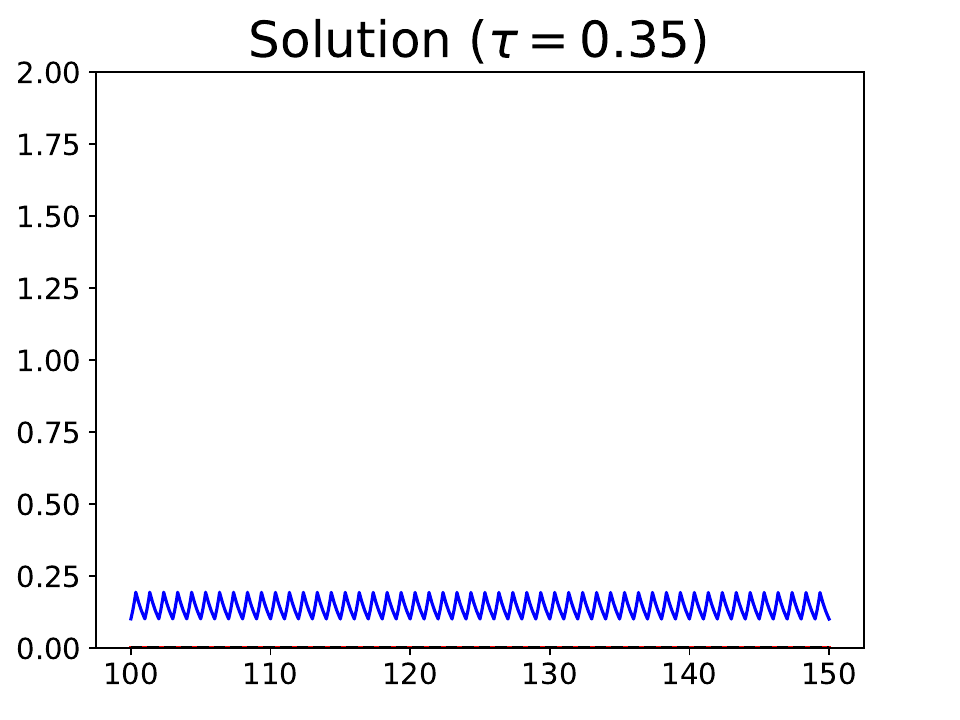}
\includegraphics[width=0.32\linewidth]{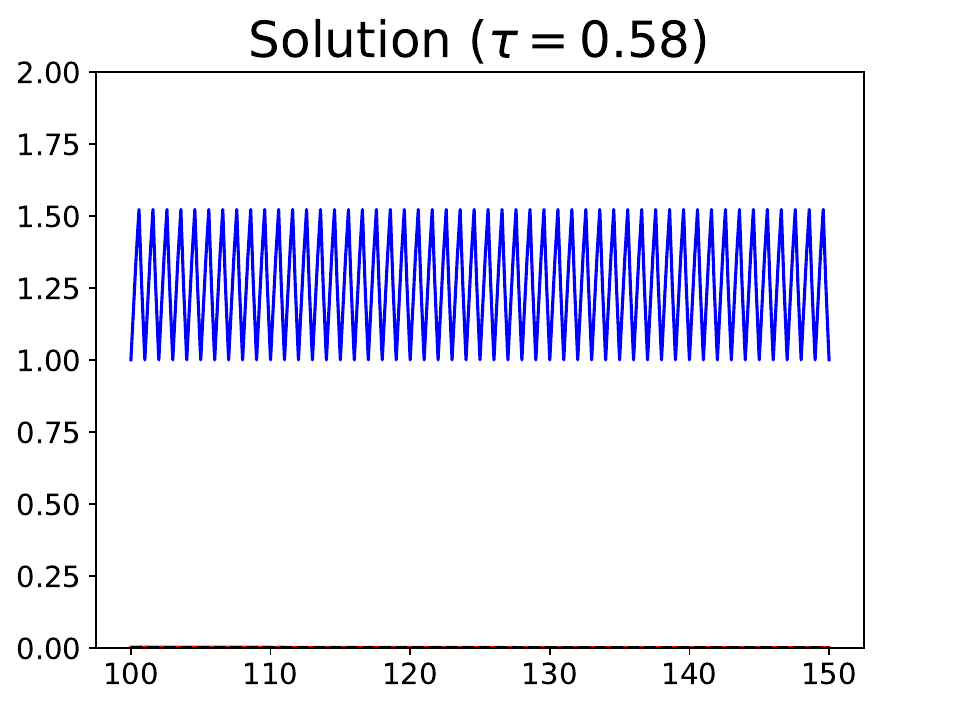}
\includegraphics[width=0.32\linewidth]{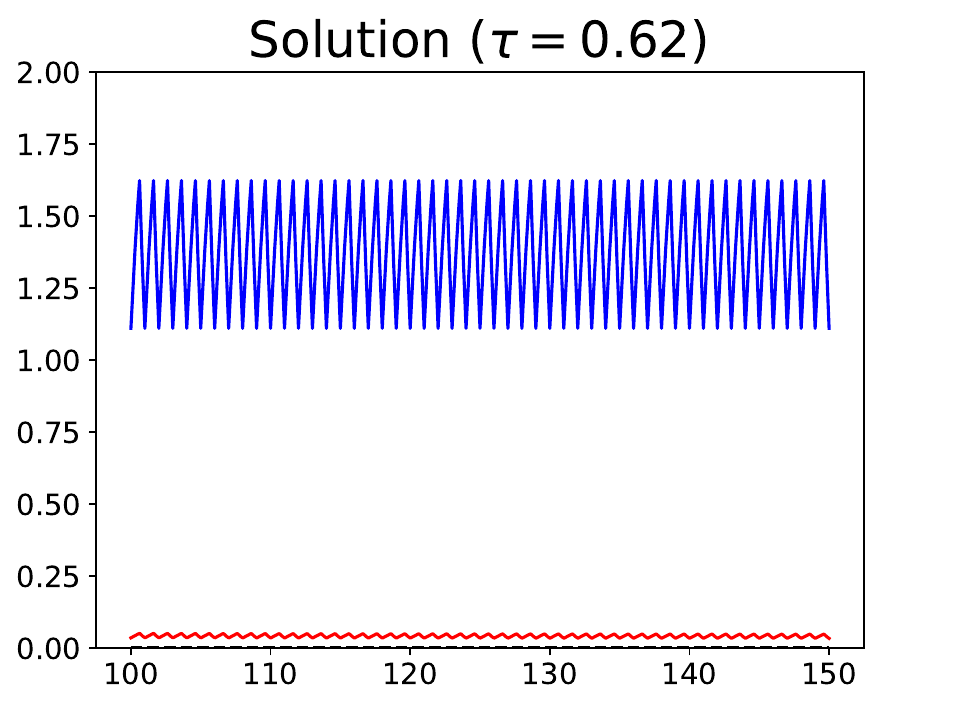}
\caption{{\small Experiment 1. Several instances of the solution of the dynamical system. Left: just after the first bifurcation.  Center: just before the second bifurcation. Right: just after the second bifurcation. The first row shows the initial periods, and the second, the final periods before reaching the equilibrium.}}
\label{fig.exp1.2}
\end{figure}

\begin{figure}[t]
\centering
\includegraphics[width=0.4\linewidth]{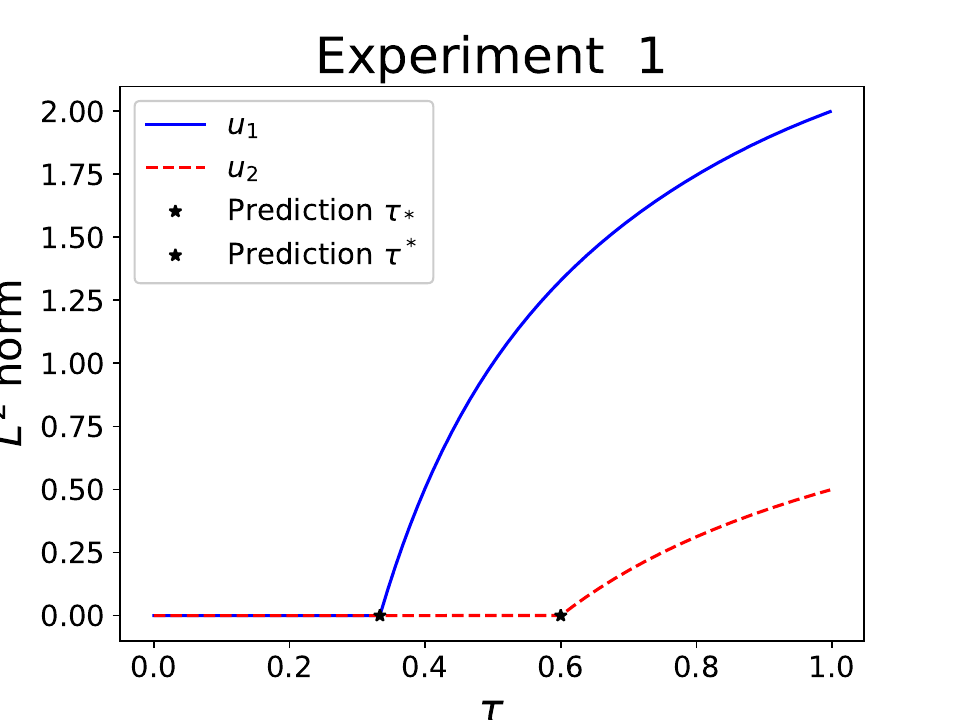}
\includegraphics[width=0.4\linewidth]{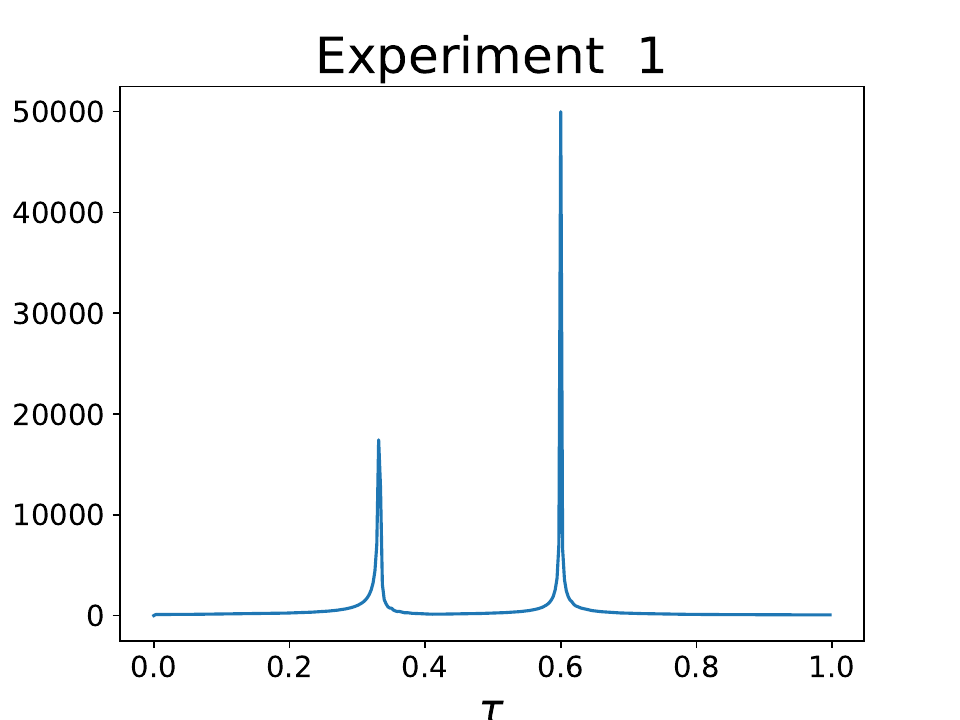}\\
\includegraphics[width=0.4\linewidth]{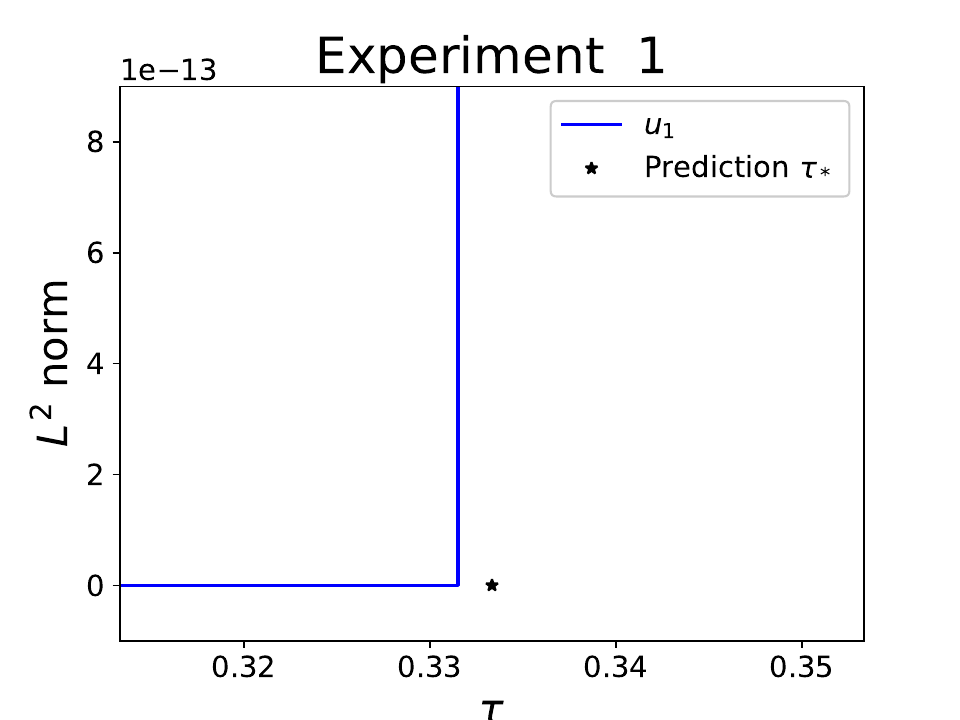}
\includegraphics[width=0.4\linewidth]{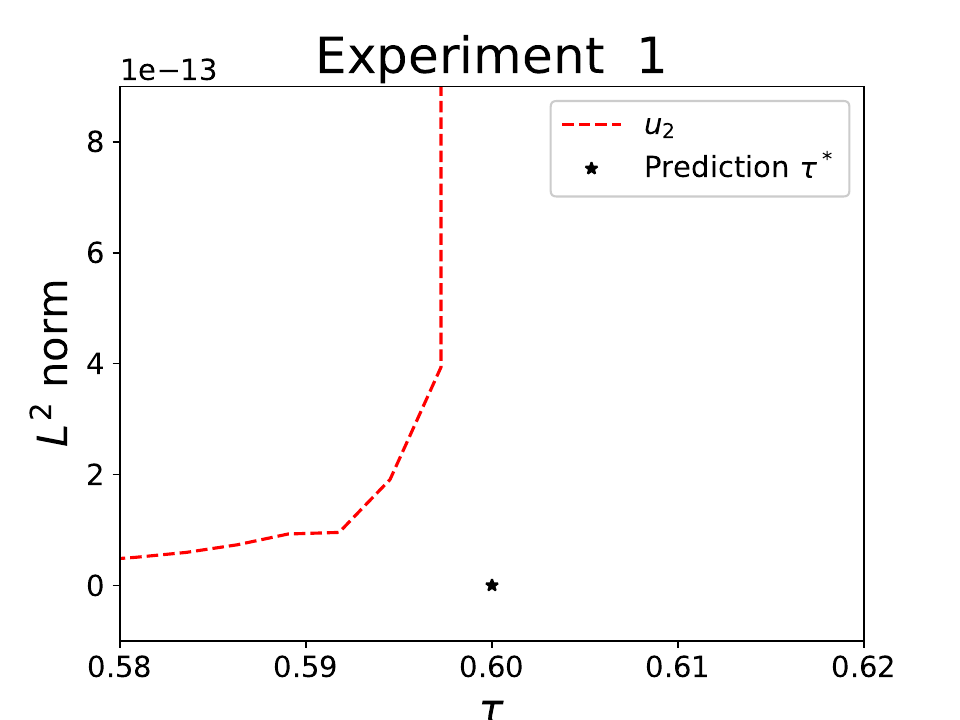}
  \caption{{\small Experiment 1: First row. Left: Bifurcation diagram, showing the expected theoretical values of $\tau_*$ and $\tau^*$, and the actual branches of computed equilibria. Right: Number of iterations to reach equilibrium as a function of $\tau$. Around both critical bifurcation parameters, convergence to equilibrium is very slow. 
  Second row: Details of the bifurcation diagram around $\tau_*$ and $\tau^*$. Notice the scale $10^{-13}$. }} 
\label{fig.exp1.3}
\end{figure}

In practice, we define the arrays \text{tau} and \texttt{U}, representing a one-period time mesh and the integral $\hat U$, respectively. Then we define 
\begin{align*}
	\texttt{tau\_star} = \argmin_{\texttt{n} \in \{0, \ldots, 365\}} \Big| \frac{\beta_{21}*\texttt{U[n]} + \mu_2}{\alpha_2+\mu_2} - \texttt{tau[n]} \Big|.
\end{align*}
In agreement with the analytical results, we obtain in both the diagonal and the non-diagonal cases the value \texttt{tau\_star = 0.6}.

\subsubsection{Experiment 1: Diagonal case $\beta_{12}=0$}

In Figure~\ref{fig.exp1.2}, we show the initial periods of the solution of the dynamical system \fer{prob2}-\fer{prob2b} and the final periods of this solution, previous to reaching the equilibrium, for three interesting season length parameters: just after the first bifurcation $\tau_*=1/3$, when the first component becomes positive. Just before the secondary bifurcation  $\tau^*=0.6$, when the second component is close to take off. And just after $\tau^*$, when the second component has taken off and becomes positive.

\begin{figure}[t]
\centering
\includegraphics[width=0.4\linewidth]{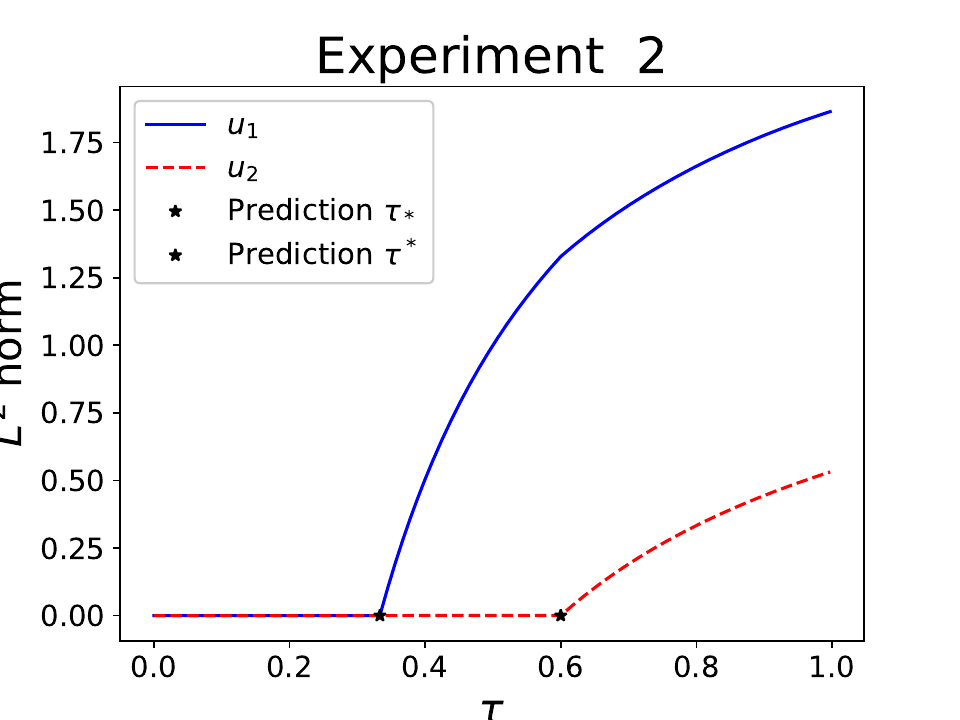}
\includegraphics[width=0.4\linewidth]{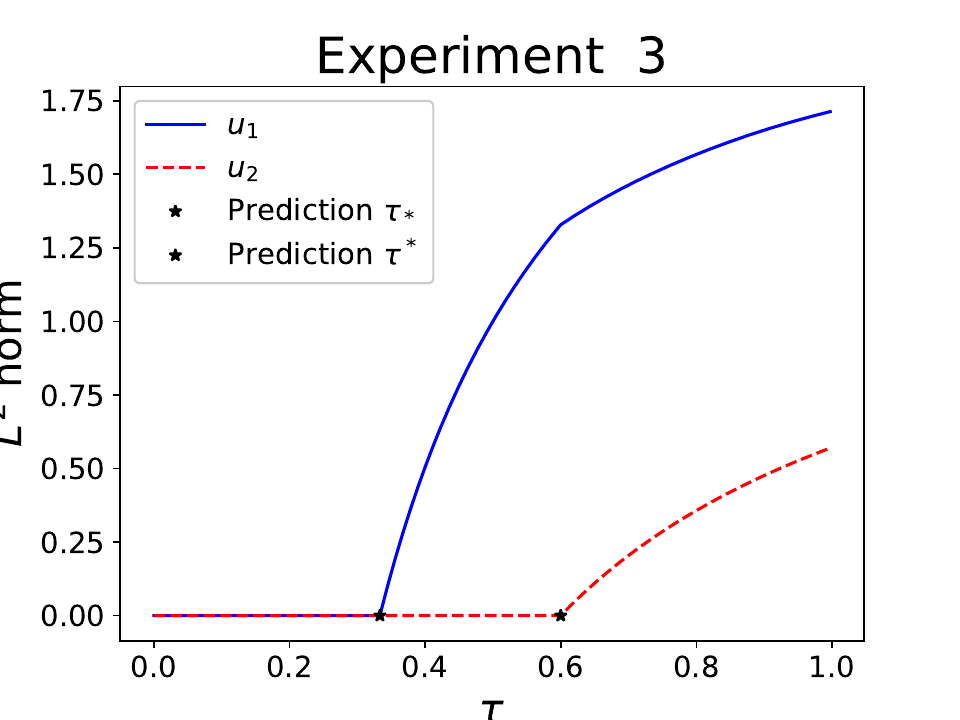}\\
\includegraphics[width=0.4\linewidth]{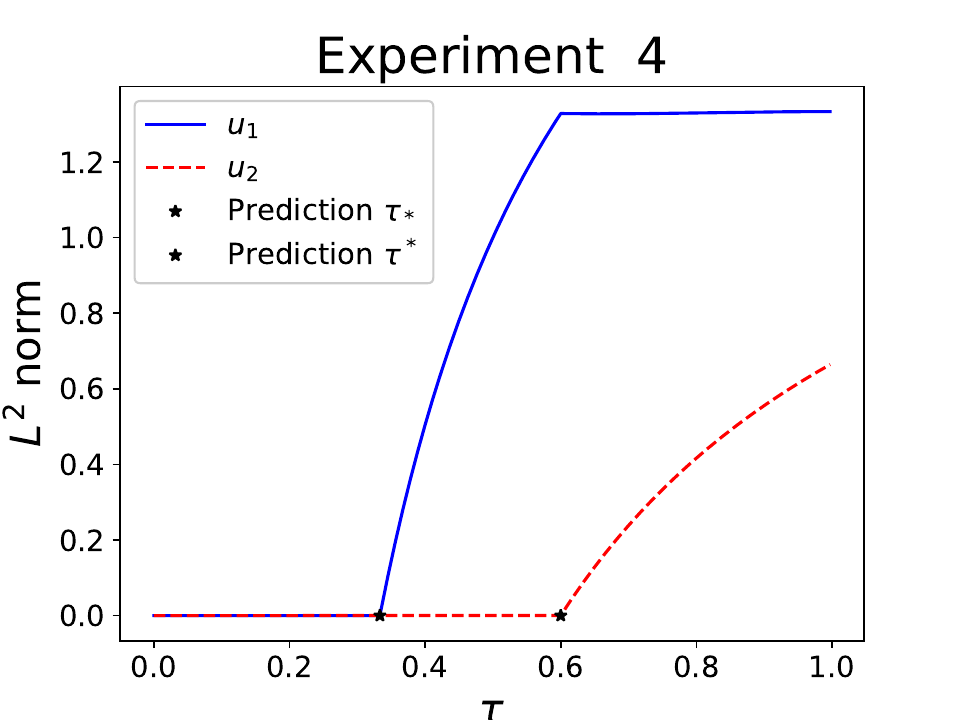}
\includegraphics[width=0.4\linewidth]{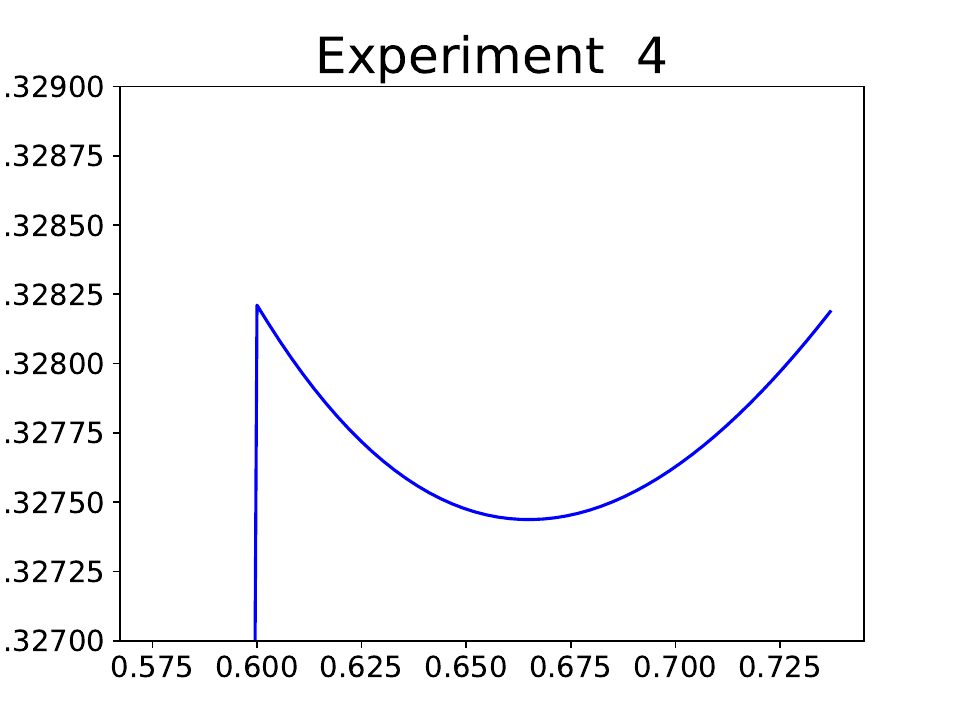}
  \caption{{\small Bifurcation diagrams for the experiments of the non-commutative case. Notice that  the component $u_1$ is influenced by the second bifurcation, see \fer{sec.bifnonsym},  developing a noticiable change of slope at $\tau^*$. The effect intensifies with increasing $\beta_{12}$, and may even lead to a reduction in the norm, as demonstrated in the bottom-right panel.}} 
\label{fig.exp2_3}
\end{figure}

In Figure~\ref{fig.exp1.3}, we show the bifurcation diagram corresponding to the $L^2(0,1)$ norm of each component of the numerically computed equilibria, as a function of the season length parameter $\tau\in(0,1)$. In the second row, we plot the details of the bifurcation diagram around the critical values $\tau_*$ and $\tau^*$, showing a reasonable agreement with the theoretical bifurcation prediction. Notice the scale $10^{-13}$.  

The number of iterations (periods) to reach the equilibrium as a function of $\tau$ is shown in Figure~\ref{fig.exp1.3}, first row, right panel.  For most  $\tau\in(0,1)$ the convergence to equilibrium is very quick. However when $\tau$ is close to both critical bifurcation parameters, the convergence speed slows down and, in fact, for $\tau\approx\tau^*$, the limit $T=T_{max}$ is reached.

\subsubsection{Experiments 2,3, and 4: Non-diagonal case $\beta_{12}>0$}

We choose the values $\beta_{12}=0.25$, $0.5$, and $1$ in the Experiments~2, 3, and 4, respectively, and conduct a similar analysis as the performed in Experiment~1.  The corresponding bifurcation diagrams are plotted in Figure~\ref{fig.exp2_3}. 

As it can be seen, the critical bifurcation parameters are again well approximated by the numerically computed equilibria. 
 
The main difference between these diagrams and that of Experiment~1 is the slope discontinuity of the bifurcation branch corresponding to $u_1$ after the secondary bifurcation $\tau^*$. 

This is in good agreement with the form of the equilibria predicted by our theoretical bifurcation analysis in which, for the diagonal case, the secondary bifurcation does not affect, locally, to $u_1$, see \fer{sec.bif}, while in the non-diagonal case, it does, see \fer{sec.bifnonsym}. As \fer{sec.bifnonsym} also predicts, this effect becomes stronger for larger values of $\beta_{12}$ (or $G_{12}[s_*](t,0)$, in that expression) and, in fact, it can reverse the growth of $\norm{u_1^*}_{L^2(0,1)}$ for a significan set of values of $\tau$.

\section{Conclusion}\label{sec.conclusion}

We have analyzed  a class of seasonally structured dynamical systems governed by a bifurcation parameter $\tau$ determining the lengths of the growth and decline seasons. Extending the existent literature, our approach allows to consider a large collection of growth and decline  functions $f_g$ and $f_d$ only restricted by some qualitative conditions on their behavior. The analysis, based on the Crandall–Rabinowitz theorem, demonstrates that a bifurcation from the extinction equilibrium occurs when a single eigenvalue of the linearized operator crosses zero and some transversality condition is fulfilled. We have further shown that secondary bifurcations may arise from the non-trivial equilibrium generated after the primary bifurcation, and we have derived sufficient conditions for their occurrence.

Numerical experiments conducted on the  paradigmatic two-species competitive Lotka-Volterra model with Malthusian decline confirm the theoretical predictions and provide precise estimates of critical season lengths. The simulations also illustrate how the simmetry (or the lack of) in interspecific interactions affects the structure of bifurcation branches, with  non-symmetric cases leading to more intricate patterns of coexistence.

Our results underline the influence of season length on the long-term behavior of populations and other seasonal systems. While primary bifurcations can be addressed analytically, secondary ones typically require numerical tools due to their dependence on unknown equilibria. 

Finally, let us emphasize that although in the example developed in Section~\ref{sec.LVM} we have taken $\varepsilon = 0$ as a limiting case representing sharp seasonal transitions, the case of $\varepsilon$ being of order one is also of interest. This corresponds to settings in which seasonal variation occurs smoothly over time, as is the case in certain epidemiological models where the infectivity rate varies gradually across the year. The analytical and numerical study of this kind of smooth seasonal transitions will be the subject of future research.


\bibliographystyle{elsarticle-num} 

\begin{thebibliography}{10}
\expandafter\ifx\csname url\endcsname\relax
  \def\url#1{\texttt{#1}}\fi
\expandafter\ifx\csname urlprefix\endcsname\relax\def\urlprefix{URL }\fi
\expandafter\ifx\csname href\endcsname\relax
  \def\href#1#2{#2} \def\path#1{#1}\fi

\bibitem{Steinthorsdottir2019}
M.~Steinthorsdottir, F.~W.-C. and, Hot summers ahead? {M}ulti-decadal spring
  season warming precedes sudden summer temperature rise in pre-anthropogenic
  climate change, GFF 141~(3) (2019) 175--180.

\bibitem{Ruosteenoja2020}
K.~Ruosteenoja, T.~Markkanen, J.~Räisänen, Thermal seasons in northern europe
  in projected future climate, Int. J. Climatol. 40~(10) (2020) 4444--4462.

\bibitem{Wang2021}
J.~Wang, Y.~Guan, L.~Wu, X.~Guan, W.~Cai, J.~Huang, W.~Dong, B.~Zhang, Changing
  lengths of the four seasons by global warming, Geophys. Res. Lett. 48~(6)
  (2021) e2020GL091753.

\bibitem{Yiqi2023}
C.~Yiqi, Z.~Yuanjie, L.~Yubin, S.~Shugang, Changes in lengths of the four
  seasons in china and the relationship with changing climate during
  1961–2020, Int. J. Climatol. 43~(3) (2023) 1349--1366.

\bibitem{Klausmeier2010}
C.~A. Klausmeier, Successional state dynamics: a novel approach to modeling
  nonequilibrium foodweb dynamics, J. Theor. Biol. 262 (2010) 584--595.

\bibitem{Hsu2012}
S.~B. Hsu, X.~Q. Zhao, A {L}otka-{V}olterra competition model with seasonal
  succession, J. Math. Biol. 64 (2012) 109--130.

\bibitem{Niu2021}
L.~Niu, Y.~Wang, X.~Xie, Carrying simplex in the {L}otka-{V}olterra competition
  model with seasonal succession with applications, Discrete Contin. Dyn.
  Syst., Ser. B 26 (2021) 2161--2172.

\bibitem{WangZhang2021}
M.~Wang, Q.~Zhang, X.~Q. Zhao, Dynamics for a diffusive competition model with
  seasonal succession and different free boundaries, J. Differ. Equ. 285 (2021)
  536--582.

\bibitem{Liu2024}
Y.~Liu, J.~Yu, J.~Li, Global dynamics of a competitive system with seasonal
  succession and different harvesting strategies, J. Differ. Equ. 382 (2024)
  211--245.

\bibitem{Pan2025}
X.~Pan, H.~Shu, L.~Wang, X.-S. Wang, J.~Yu, On the periodic solutions of
  switching scalar dynamical systems, J. Differ. Equ. 415 (2025) 365--382.

\bibitem{Namba1986}
T.~Namba, Bifurcation phenomena appearing in the {L}otka-{V}olterra competition
  equations: a numerical study, Math. Biosci. 81~(2) (1986) 191--212.

\bibitem{Peng2012}
R.~Peng, X.~Q. Zhao, A reaction-diffusion sis epidemic model in a time-periodic
  environment, Nonlinearity 25 (2012) 1451--1471.

\bibitem{Lou2017}
Y.~Lou, X.~Q. Zhao, A theoretical approach to understanding population dynamics
  with seasonal developmental durations, J. Nonlinear Sci. 27 (2017) 573--603.

\bibitem{Browne2020}
C.~J. Browne, X.~Pan, H.~Shu, X.~Wang, Resonance of periodic combination
  antiviral therapy and intracellular delays in virus model, Bull. Math. Biol.
  82 (2020) 1--29.

\bibitem{Chen2023}
J.~Chen, J.~O.~R. Rincon, G.~DeGrandi-Hoffman, J.~Fewell, J.~Harrison, Y.~Kang,
  Impacts of seasonality and parasitism on honey bee population dynamics, J.
  Math. Biol. 87 (2023) 19.

\bibitem{Seeley1985}
T.~D. Seeley, P.~K. Visscher, Survival of honeybees in cold climates: the
  critical timing of colony growth and reproduction, Ecol. Entomol. 10~(1)
  (1985) 81--88.

\bibitem{Almeida2018}
A.~Almeida, C.~Codeço, P.~M. Luz, Seasonal dynamics of influenza in brazil:
  the latitude effect, BMC Infect. Dis. 18~(695) (2018).

\bibitem{Lou2019}
Y.~Lou, K.~Liu, D.~He, D.~Gau, S.~Ruan, Modelling diapause in mosquito
  population growth, J. Math. Biol. 28 (2019) 2259--2288.

\bibitem{Bai2020}
Z.~Bai, Y.~Lou, X.-Q. Zhao, A delayed succession model with diffusion for the
  impact of diapause on population growth, SIAM J. Appl. Math. 80 (2020)
  1493--1519.

\bibitem{Namba1984}
T.~Namba, Competitive co-existence in a seasonally fluctuating environment, J.
  Theor. Biol. 111~(2) (1984) 369--386.

\bibitem{Namba1993}
T.~Namba, S.~Takahashi, Competitive coexistence in a seasonally fluctuating
  environment {II}. multiple stable states and invasion success, Theor. Pop.
  Biol. 44 (1993) 374--402.

\bibitem{Hu1995}
S.~S. Hu, A.~J. Tessier, Seasonal succession and the strength of intra- and
  interspecific competition in a {D}aphnia assemblage, Ecology 76 (1995)
  2278--2294.

\bibitem{DuBowy1988}
P.~J. DuBowy, Waterfowl communities and seasonal environments: temporal
  variability in interspecific competition, Ecology 69 (1988) 1439--1453.

\bibitem{DeAngelis2009}
D.~L.~D. Angelis, J.~C. Trexler, D.~D. Donalson, Competition dynamics in a
  seasonally varying wetland, in: S.~Cantrell, C.~Cosener, S.~Ruan (Eds.),
  Spatial ecology, CRC Press/Chapman and Hall, London, 2009, pp. 1--13.

\bibitem{Koch1974}
A.~L. Koch, Coexistence resulting from an alternation of density dependent and
  density independent growth, J. Theor. Biol. 44 (1974) 373--386.

\bibitem{Liu2007}
P.~Liu, J.~Shi, Y.~Wang, Imperfect transcritical and pitchfork bifurcations, J.
  Funct. Anal. 251~(2) (2007) 573--600.

\bibitem{Crandall1971}
M.~G. Crandall, P.~H. Rabinowitz, Bifurcation from simple eigenvalues, J.
  Funct. Anal. 8~(2) (1971) 321--340.

\bibitem{Crandall1973}
M.~G. Crandall, P.~H. Rabinowitz, Bifurcation, perturbation of simple
  eigenvalues, and linearized stability, Arch. Ration. Mec. An. 52 (1973)
  161--180.

\bibitem{Liu2013}
P.~Liu, J.~Shi, Y.~Wang, Bifurcation from a degenerate simple eigenvalue, J.
  Funct. Anal. 264~(10) (2013) 2269--2299.

\bibitem{Liu2013b}
P.~Liu, J.~Shi, Y.~Wang, A double saddle-node bifurcation theorem, Commun. Pure
  Appl. Anal. 12~(6) (2013) 2923--2933.

\bibitem{Brezisbook}
H.~Brezis, Functional Analysis, Sobolev Spaces and Partial Differential
  Equations, Springer, 2011.

\bibitem{Scipy2020}
P.~Virtanen, R.~Gommers, T.~E. Oliphant, M.~Haberland, T.~Reddy, D.~Cournapeau,
  E.~Burovski, P.~Peterson, W.~Weckesser, J.~Bright, S.~J. {van der Walt},
  M.~Brett, J.~Wilson, K.~J. Millman, N.~Mayorov, A.~R.~J. Nelson, E.~Jones,
  R.~Kern, E.~Larson, C.~J. Carey, {\.I}.~Polat, Y.~Feng, E.~W. Moore,
  J.~{VanderPlas}, D.~Laxalde, J.~Perktold, R.~Cimrman, I.~Henriksen, E.~A.
  Quintero, C.~R. Harris, A.~M. Archibald, A.~H. Ribeiro, F.~Pedregosa, P.~{van
  Mulbregt}, {SciPy 1.0 Contributors}, {{SciPy} 1.0: Fundamental Algorithms for
  Scientific Computing in Python}, Nat. Methods 17 (2020) 261--272.

\end{thebibliography}


\end{document}